[1994/12/01]
\documentclass[10pt, reqno]{amsart}
\usepackage{a4wide}
\usepackage[english, activeacute]{babel}
\usepackage{amsmath,amsthm,amsxtra}
\usepackage{amssymb}
\usepackage{latexsym}
\usepackage{amsfonts}
\usepackage{hyperref}
\usepackage{pdftricks}

\begin{psinputs}
\usepackage{pst-plot}
\end{psinputs}

\title{Well-posedness and stability results for the Gardner equation}
\author{Miguel A. Alejo}
\address{Departament of Mathematics, University of the Basque Country, Bilbao \\ Spain}
\email{miguelangel.alejo@ehu.es}
\date{December, 2010}
\thanks{}

\chardef\bslash=`\\ 

\newtheorem{thm}{Theorem}[section]

\newtheorem{lem}[thm]{Lemma}

\theoremstyle{definition}

\theoremstyle{remark}

\numberwithin{equation}{section}

\newcommand{\bel}{\begin{equation}\label}
\newcommand{\eeq}{\end{equation}}
\newcommand{\R}{\mathbb{R}}

\newcommand{\T}{\mathbb{T}}

\newcommand{\C}{\mathbb{C}}
\newcommand{\lm}{\lesssim}

\newcommand{\supp}{\operatorname{supp}}

\def\bm{\left( \begin{array}{cc}}
\def\endm{\end{array}\right)}

\newcommand{\be}{\begin{equation}}
\newcommand{\ee}{\end{equation}}
\newcommand{\ba}{\begin{equation*}}
\newcommand{\ea}{\begin{equation*}}
\newcommand{\bea}{\begin{eqnarray}}
\newcommand{\eea}{\end{eqnarray}}
\newcommand{\bee}{\begin{eqnarray*}}
\newcommand{\eee}{\end{eqnarray*}}
\newcommand{\ben}{\begin{enumerate}}
\newcommand{\een}{\end{enumerate}}

\newcommand{\eval}[2][\right]{\relax
  \ifx#1\right\relax \left.\fi#2#1\rvert}

\begin{document}
\begin{abstract}
In this article we present local well-posedness results in the classical Sobolev space $H^s(\R)$ with $s>1/4$
for the Cauchy problem of the Gardner equation, overcoming the problem of the loss of the scaling property of this equation. We also cover the energy space $H^1(\R)$ where global well-posedness follows from the conservation laws of the system. Moreover, we construct  solitons of the Gardner equation explicitly  and  prove that, under certain conditions,  this family is orbitally stable in the energy space.
\end{abstract}
\maketitle \markboth{Well-posedness and stability results for the Gardner equation} {Miguel A. Alejo}
\renewcommand{\sectionmark}[1]{}

\section{{\bf Introduction}}\label{intro}

\medskip

In this article we present some results about the regularity  for the Cauchy problem of the  Gardner equation (shortly GE) 

\bel{GE0}
\begin{cases}
v_t + v_{xxx} + 6\sigma(v^2)_x+2(v^3)_x=0~~~~~~\sigma, t, x\in\mathbb{R},\\
v(x,0)=v_0(x)\in H^s(\mathbb{R}),
\end{cases}\eeq
with data in the classical Sobolev space $H^s(\R)$ and the orbital stability of its soliton solutions.  In a previous work \cite{Ale} we faced the problem of finding some $L^{\infty}$-solutions with nonzero limits at infinity and with geometrical interpretation of the focusing {\it modified} Korteweg-de Vries (shortly mKdV) 
\bel{mKdVf}u_t + u_{xxx} + 2(u^3)_x=0,\eeq
that is, solutions of the type $\sigma+v(x+ct)$, with $c>0$, $\sigma\in\R$ and  $v$ a traveling wave solution with exponential decay at infinity. When  introducing such ansatz in the focusing mKdV \eqref{mKdVf}, the GE \eqref{GE0}  appears (up to rescaling). This evolution equation is characterized to be composed by a KdV term $(v^2)_x$ and a positive mKdV term $(v^3)_x$. The competition between these two different nonlinear terms together with the linear  dispersive term $v_{xxx}$ allows the existence of more intricate soliton solutions (see section \ref{Stability}) as well as exact breather solutions (see \cite{Ale}). The Gardner equation plays also an important role in the proof of the $L^2-$stability of the multisoliton solution of KdV, through to the so called Gardner transform which links $H^1-$solutions  of the Gardner equation with $L^2-$solutions of the KdV equation (see \cite{AMV}).\\

Soliton solutions (or equivalently $\sigma$-soliton solutions of the mKdV, if their asymptotic constant is equal to $\sigma$) in the focusing and defocusing cases ($\pm$ sign respectively in the cubic nonlinearity of \eqref{mKdVf}) are easily related through the transformation of the asymptotic parameter $\sigma$ to $i\sigma$. Indeed, these solitons can be explicitly obtained integrating the resulting second order ODE which arises when we look for traveling wave solutions of \eqref{GE0} of the type $v(x+ct)$ with $c>0$ (see section \ref{Stability}).\\

In this paper we are also interested to prove the orbital stability {\it \'a la Zhidkov} of these soliton solutions under small perturbations in $H^1(\R)$. Hence, we need a global well posedness (GWP for short) result for the initial value problem (IVP for short) \eqref{GE0} in the energy space $H^1(\R)$, and therefore it is enough for our aim to prove the local well posedness (LWP for short) of the IVP \eqref{GE0} below the conservation law $H^1(\R)$. Indeed, we prove the LWP in $H^s(\R)$ with $s>1/4$.  Note that the Gardner equation is not scaling invariant. This loss of the scale property raises two problems. The first one appears in the proof of persistence of the solution in the LWP result since the scale property can  not be used and we need introduce a rescaling of the problem in terms of a new auxiliar function. The second problem  is that the proof of the convexity condition suggested by P. Zhidkov can not be deduced directly and we need to integrate the Lyapunov functional. As a consequence, we see that the proof of the stability of solitons in the focusing and defocusing cases are almost identical. Therefore, we will only show here the focusing case. \\

The LWP of the IVP for KdV with initial data in $H^{s}(\R),s>-3/4$ was obtained by  C.Kenig, G.Ponce and L.Vega  in \cite{KePV4}. They showed sharp bilinear estimates in the functional space $X^{s,b}$, introduced by Bourgain in \cite{B}, up to the index  $s=-3/4$. In  \cite{CtCoTa}, M. Christ, J. Colliander and T. Tao proved the LWP of the IVP for KdV with initial data in $H^s(\R),~~s\geq-3/4$, using a {\it modified} Miura transform and the existence theory for the mKdV. They also proved the global theory for initial data in  $H^s(\R),~~s>-3/4$.\\ 

The LWP of the IVP for the focusing mKdV with inital data in the Sobolev space $H^{s}(\R),~s\geq1/4$, was given by C.Kenig, G.Ponce and L.Vega in \cite{KePV2}, where they also proved the global well posedness in the energy space $H^s(\R),~s\geq1$. The global result below the conservation law was shown, for initial data in $H^s(\R),~~s>1/4$  by J. Colliander, M. Keel, G. Staffilani, H. Takaoka and T. Tao in  \cite{CoKeTakTa}, using the existence theory for KdV and the Miura transform.\\

In the next sections we present the local well-posedness of the IVP for the GE \eqref{GE0} in the Sobolev space $H^s(\R),s>1/4$, obtaining bilinear and trilinear estimates by using the auxiliar space of functions  $X^{s,b}(\R\times\R)$ and the $[k;\R]-$multiplier theory introduced by T.Tao \cite{T}. We also cover the energy space $H^1(\R)$ where global well posedness follows from the conservation laws of the problem. Moreover, we construct solitons of the GE and we prove, under certain conditions, that this family is orbitally stable in $H^1(\R)$.

\section{\textbf{Local Theory}}\label{localWP}
We are interested in the local (in time, $t=T>0$) existence theory for the IVP of the GE \eqref{GE0} with initial data in\footnote[1]{In what follows, $H^{1/4^{+}}(\R)$ means $H^s(\R)$ with $s>1/4$.} $H^{1/4^{+}}(\R)$. Since the GE is not invariant under scaling transformations, we can not proceed as usual, i.e., prove the LWP of \eqref{GE0} at time $t=1$ and use the scaling transformation to extent the local time to $t=T>0$. So we devise a procedure closely related to the usual one, introducing an auxiliar function $w$ through a scaling transformation. Next, we  prove the LWP at $t=1$ of the associated GE for $w$ at a certain scaling parameter $\lambda$ (to be determined). Finally we use this result and the scaling relation to obtain directly the LWP of the IVP for the GE \eqref{GE0} at $t=T>0$.\\

We make the following  scaling transformation: 
\bel{escala1}v(x,t)=\lambda^\alpha w(\lambda x,\lambda³t),~\lambda>0,~\alpha>2,\eeq
where $v$ is a solution of the IVP for the GE \eqref{GE0}, so that the associated Gardner equation for the auxiliar function $w$ is

\bel{GEw}
\begin{cases}
w_t + w_{xxx} + 6\sigma\lambda^{\alpha-2}(w^2)_x+2\lambda^{2(\alpha-1)}(w^3)_x=0,~~~~~~\sigma,
\lambda>0 , \alpha,t, x\in\mathbb{R},\\
w(x,0)=w_0(x)\in H^s(\mathbb{R}).
\end{cases}\eeq

For $s,b\in\R$, the space of functions $X^{s,b}$ denotes the completion of the Schwartz space $\mathcal{S}(\R^2)$ with respect to the norm
\bel{normaXsb}||f||_{X^{s,b}}=\left(\int_{-\infty}^{+\infty}\int_{-\infty}^{+\infty}(1+|\tau-\xi^3|)^{2b}(1+|\xi|)^{2s}|\hat{f}(\xi,\tau)|^2d\xi d\tau\right)^{1/2}.\eeq
We invoke the work of T.Tao \cite{T} about multilinear estimates to define a $[k;\R]$-multiplier as any function $m:\Gamma_k(\R)\rightarrow\C$, where $\Gamma_k(\R)$ is the hyperplane

\bel{desiTao1}\Gamma_k(\R):=\{(\xi_1,\xi_2,\dots,\xi_k)\in \R^k : \xi_1+\xi_2+\dots+\xi_k=0\},\eeq
endowed with the measure

\bel{desiTao2}\int_{\Gamma_k(\R)}f:=\int_{\R^{k-1}}f(\xi_1,\dots,\xi_{k-1},-\xi_1-\dots-\xi_k-1)d\xi_1d\xi_2\dots d\xi_{k-1}.\eeq

For $||m||_{[k;\R]}$ we denote the best constant such that the inequality

\bel{desiTao}\left|\int_{\Gamma_k(\R)}m(\xi)\prod^k_{j=1}f_j(\xi_j)\right|\leq||m||_{[k; \R]}\prod^k_{j=1}||f_j||_{L^2(\R)},\eeq
holds for all test functions $f_i$ on $\R$. In the sequel, the free operator $W(t)$ denotes an element of the unitary group $\{W(t)\}^{+\infty}_{-\infty}$ describing the solution of the linear IVP associated to \eqref{GE0}

\bel{linealIVP}
\begin{cases}
v_t + v_{xxx}=0,~~~t, x\in\mathbb{R},\\
v(x,0)=v_0(x),
\end{cases}\eeq
where 
\bel{solLinearIVP}v(x,t)=W(t)v_0(x)=\int_{\R}e^{ix\xi+t\xi^3}\hat{v}_0(\xi)d\xi.\eeq
Finally, let $\psi\in C^{\infty}_0(\R)$ with $\psi\equiv1$ on $[-1/2,1/2]$ and $\supp\psi\subseteq(-1,1)$. Then, our main results are:

\begin{thm}\label{LWPfw} Let be $s>1/4$ and  $\sigma>0$. Then there exist constants $d_1>0,~~d_2>0,~~\lambda>0$ and
$b\in (1/2, 1)$, such that for all $w_0\in H^{s}(\R)$, exists a unique solution $w \in C([-1, 1]:H^s(\R))$ of \eqref{GEw}
%
with $\lambda\leq\min(d_1||w_0||^{\frac{-1}{\alpha-2}}_{H^s},d_2||w_0||^{\frac{-1}{\alpha-1}} _{ H^s}),$ and $w$  satisfying

\bel{persistencia1w}w\in C([-1,1]:H^s(\mathbb{R})),\eeq

\bel{persistencia2w}w\in X^{s,b}\subseteq
L^{p}_{x,~loc}(\mathbb{R}:L^2_t([-1,1])),~1\leq p\leq\infty,\eeq

\bel{persistencia3w}\partial_x(w^3)\in X^{s,b-1}(\R\times\R)~~ \wedge~~
\partial_x(w^2)\in X^{s,b-1}(\R\times\R).\eeq
\end{thm}

\begin{thm}\label{LWPf} Let be $s>1/4$ and $\sigma>0$. Then there exists $b\in(1/2,1)$, such that for every  $v_0(x)\in H^s(\mathbb{R})$, there exist a local time  $T=T(||v_0||_{H^s})>0$ (with
$T(\rho)\rightarrow\infty$ when $\rho\rightarrow0$) and a unique solution $v(x,t)\equiv v(t)$ of \eqref{GE0} such that $\sigma + v(t)$ is the unique solution of \eqref{mKdVf} in the interval $[-T,T]$ satisfying\\
\bel{persistencia1}v\in C([-T,T]:H^s(\mathbb{R})),\eeq

\bel{persistencia2}v\in X^{s,b}\subseteq
L^{p}_{x,~loc}(\mathbb{R}:L^2_t([-T,T])),~1\leq p\leq\infty,\eeq

\bel{persistencia3}\partial_x(v^3)\in X^{s,b-1}(\R\times\R)~~\wedge~~
\partial_x(v^2)\in X^{s,b-1}(\R\times\R).\eeq

\end{thm}


We will resort to the following preliminary estimates for the free operator $W(t)$, of the unitary group describing the solution of the linear IVP associated to \eqref{GE0} (for proofs of such estimates, see \cite{KePV3}):\\

\begin{lem}\label{lemaa} Let $s>1/4,~1/2<b<1$ and $0<\delta<1$. Then, there exists a constant $c>0$ such that:
\begin{enumerate}
 \item \label{lema1}$||\psi(\delta^{-1}t)W(t)v_0||_{X^{s,b}}\leq
c~\delta^{\frac{1-2b}{2}}||v_0||_{H^{s}}.$\\
\item \label{lema2}$||\psi(\delta^{-1}t)h||_{X^{s,b}}\leq
c~\delta^{\frac{1-2b}{2}}||h||_{X^{s,b}}.$\\
\item \label{lema3}
$||\psi(\delta^{-1}t)\int^t_0W(t-t')w(t')dt'||_{X^{s,b}}\leq
c~\delta^{\frac{1-2b}{2}}||w||_{X^{s,b-1}}.$\\
\item \label{lema4}
$||\psi(\delta^{-1}t)\int^t_0W(t-t')w(t')dt'||_{H^{s}}\leq
c~\delta^{\frac{1-2b}{2}}||w||_{X^{s,b-1}}.$
\end{enumerate}
\end{lem}






Moreover, to prove the local theorem we will need the following bilinear and trilinear estimates:
\begin{lem}\label{lemabilinear}
Let $s>1/4$. Then, for all $u_i=\psi(t)\phi_i(x,t),~i=1,2,$ with support in $\R\times[-1,1]$ and $b=1/2+\epsilon,~~0<\epsilon\ll1$, the following inequality holds:
\bel{bilinear} ||u_1u_2||_{L^2(\R\times\R)}\leq c~ ||\phi_1||_{X^{s,b}(\R\times\R)} ||\phi_2||_{X^{-1/2,1-b}(\R\times\R)}.\eeq
\end{lem}
{\it Proof.} For \eqref{bilinear} see Proposition 6.2. in \cite{T}. \begin{flushright}$\Box$\end{flushright}
\begin{lem}\label{lematrilinear}
Let $s>1/4$ and $b=1/2+\epsilon,~~0<\epsilon\ll1$. Then, for all $u_i=\psi(t)\phi_i(x,t),~i=1,2,3,$ with support in $\R\times[-1,1]$, and $c>0$,
the following inequality holds

\bel{trilinear} ||\partial_x(u_1u_2u_3)||_{X^{s,b-1}(\R\times\R)}\leq
c~||\phi_1||_{X^{s,b}(\R\times\R)}
||\phi_2||_{X^{s,b}(\R\times\R)} ||\phi_3||_{X^{s,b}(\R\times\R)}.\eeq
\end{lem}
{\it Proof\footnote{In this proof we follow the notation of \cite{T}. Here $\hat{\cdot}$ denotes the Fourier transform.}.} We emphasize that  essentially the proof appears in the work of T.Tao \cite{T}(it is worth to note that it also works for $s=1/4$), but we present here for the sake of completeness. We  summarize the proof in three steps:

 1.- By the duality of the
spaces $X^{s,b-1}(\R\times\R)$,  $X^{-s,1-b}(\R\times\R)$ and using Plancherel, we obtain that

\bel{convo1}\begin{array}{ll}\int_{\R}\int_{\R}\bar{f}
\partial_x(u_1u_2u_3)dxdt=\int_{\R}\int_{\R}\hat{\bar{f}}(-\xi,-\tau)\widehat{
\partial_x(u_1u_2u_3)}(\xi,\tau)d\xi d\tau\\\\
=\int_{\R}\int_{\R}(i\xi)\hat{\bar{f}}(-\xi,-\tau)\widehat{u_1u_2u_3}(\xi,
\tau)d\xi d\tau.\end{array}\eeq
Then, taking into account that $\xi_3=\xi-\xi_1-\xi_2,~\tau_3=\tau-\tau_1-\tau_2$, the expression \eqref{convo1} reduces to

\bel{convo6}\begin{array}{l}\int_{\R}\int_{\R}\bar{f}
\partial_x(u_1u_2u_3)dxdt=\int_{\Gamma_{3+1}(\R\times\R)}
i(\xi_1+\xi_2+\xi_3)\cdot\left(\prod_{j=1}^{3}\hat{u}_j(\xi_j,
\tau_j)\right)\cdot \hat{\bar{f}}(\xi_4,\tau_4).\end{array}\eeq
Resorting to Lemma \ref{lemaa} \eqref{lema2} and in terms of the $X^{s,b}$ norm, we have 

\bel{convo7}\begin{array}{lll}\left|\int_{\Gamma_{3+1}(\R\times\R)}
(\xi_1+\xi_2+\xi_3)\cdot\left(\prod_{j=1}^{3}\hat{u}_j(\xi_j,\tau_j)\right)\cdot
\hat{\bar{f}}(\xi_4,\tau_4)\right|\\\\
\leq c~\left\|\frac{(\xi_1+\xi_2+\xi_3)\left\langle
\xi_4\right\rangle^s\left\langle
\tau_4-\xi_4^3\right\rangle^{b-1}}{\prod^3_{j=1}\left\langle
\xi_j\right\rangle^s\left\langle \tau_j-\xi_j^3\right\rangle^b}\right\|_{[3+1;
\R\times\R]}\cdot\left(\prod^3_{j=1}||\phi_j||_{X^{s,b}}\right)||f||_{X^{-s,1-b}
}.\end{array} \eeq
So, denoting $b=1/2+\epsilon,$ for $0<\epsilon\ll1$, it is enough to prove that 
 
\bel{estimacióntrilineal2}
\left\|\frac{(\xi_1+\xi_2+\xi_3)\left\langle \xi_4\right\rangle^s\left\langle
\tau_4-\xi_4^3\right\rangle^{b-1}}{\prod^3_{j=1}\left\langle
\xi_j\right\rangle^s\left\langle\tau_j-\xi_j^3\right\rangle^{b}}\right\|_{[3+1; \R\times\R]}\lm1. \eeq\\

2.- Since $\xi_4=-\xi_1-\xi_2-\xi_3$ implies $|\xi_1+\xi_2+\xi_3|\sim\xi_4$, and applying the inequality

\bel{desifraccionaltri}\left\langle \xi_4\right\rangle^{s+1}\lm\left\langle
\xi_4\right\rangle^{1/2}\sum^3_{j=1}\left\langle \xi_j\right\rangle^{s+1/2},\eeq
\eqref{estimacióntrilineal2} simplifies as follows:

\bel{estimacióntrilineal21}\begin{array}{ll}
\left\|\frac{(\xi_1+\xi_2+\xi_3)\left\langle \xi_4\right\rangle^s\left\langle
\tau_4-\xi_4^3\right\rangle^{b-1}}{\prod^3_{j=1}\left\langle
\xi_j\right\rangle^s\left\langle
\tau_j-\xi_j^3\right\rangle^{b}}\right\|_{[3+1; \R\times\R]}\sim
\left\|\frac{\left\langle \xi_4\right\rangle^{s+1}\left\langle
\tau_4-\xi_4^3\right\rangle^{b-1}}{\prod^3_{j=1}\left\langle
\xi_j\right\rangle^s\left\langle
\tau_j-\xi_j^3\right\rangle^{b}}\right\|_{[3+1; \R\times\R]}\\\\

\lm\left\|\frac{\left\langle \xi_4\right\rangle^{1/2}\sum^3_{i=1}\left\langle
\xi_i\right\rangle^{s+1/2}\left\langle
\tau_4-\xi_4^3\right\rangle^{b-1}}{\prod^3_{j=1}\left\langle
\xi_j\right\rangle^s\left\langle
\tau_j-\xi_j^3\right\rangle^{b}}\right\|_{[3+1; \R\times\R]}.\end{array}\eeq\\

3.- Assuming now (w.l.g.) that dual variable $\xi_2$ is the greater one, we use the following estimate:

$$\frac{\sum^3_{j=1}\left\langle
\xi_j\right\rangle^{s+1/2}}{\left\langle \xi_2\right\rangle^{s}}\lm\left\langle
\xi_2\right\rangle^{1/2}.$$

Then \eqref{estimacióntrilineal21} remains as

\bel{estimacióntrilineal3}\left\|\frac{\left\langle
\xi_4\right\rangle^{1/2}\left\langle \xi_2\right\rangle^{1/2}\left\langle
\xi_1\right\rangle^{-s}\left\langle \xi_3\right\rangle^{-s}\left\langle
\tau_4-\xi_4^3\right\rangle^{b-1}}{\left\langle
\tau_2-\xi_2^3\right\rangle^{1-b}\prod^2_{i=1}\left\langle
\tau_{2i-1}-\xi_{2i-1}^3\right\rangle^{b}}\right\|_{[3+1; \R\times\R]}.
\eeq
Selecting $m$ as
$$m(\xi_1,\xi_2)=\frac{\left\langle \xi_1\right\rangle^{-s}\left\langle
\xi_2\right\rangle^{1/2}}{\left\langle
\tau_2-\xi_2^3\right\rangle^{1-b}\left\langle
\tau_1-\xi_1^3\right\rangle^{b}},$$
it is possible to rewrite \eqref{estimacióntrilineal3} as
$$\left\|m(\xi_1,\xi_2)\cdot\overline{m(-\xi_3,-\xi_4)}
\right\|_{[3+1; \R\times\R]}=\left\|m(\xi_1,\xi_2)\right\|^2_{[2+1; \R\times\R]},$$
where we have used the estimate $TT^{*}$ given in \cite[p.8]{T}. In concluding, we need that
\bel{estimacióntrilineal4}\left\|m(\xi_1,\xi_2)\right\|_{[2+1;
\R\times\R]}=\left\|\frac{\left\langle \xi_1\right\rangle^{-s}\left\langle
\xi_2\right\rangle^{1/2}}{\left\langle
\tau_2-\xi_2^3\right\rangle^{1-b}\left\langle
\tau_1-\xi_1^3\right\rangle^{b}}\right\|_{[2+1; \R\times\R]}\lm1,
\eeq
which is proved in Lemma \ref{lemabilinear}. \begin{flushright}$\Box$\end{flushright}
By using the same steps than in the proof of  \eqref{trilinear}, it is straighforward to prove the following:
\begin{lem}\label{bilinearkdv}
Let $s>0$ and $b=1/2+\epsilon,~~0<\epsilon\ll1$. Then for all
$u_i=\psi(t)\phi_i(x,t),~i=1,2,$ with support in $\R\times[-1,1]$, and $c>0$, the following inequality holds:
\bel{bilinearkdv0} ||\partial_x(u_1u_2)||_{X^{s,b-1}(\R\times\R)}\leq c~||\phi_1||_{X^{s,b}(\R\times\R)} ||\phi_2||_{X^{s,b}(\R\times\R)}.\eeq\end{lem}
%
%
{\bf Proof of Theorem \ref{LWPfw}.}\\
We denote the $H^s(\R)$ norm of the initial data $w_0(x)$ as
\bel{normainicial}~||w_0||_{H^{s}}=r_w.\eeq
 For $w_0\in H^{s}(\R)$ with $s>1/4$,  the localized Duhamel operator  is
\begin{align}\nonumber \Phi_{1,w_0}(w)\equiv\Phi_1(w)=&\psi(t)W(t)w_0 -6\sigma\lambda^{\alpha-2}\psi(t)\int^t_0W(t-t')\partial_x[(\psi(t')w(t'))^2]
dt'\\\nonumber
&-2\lambda^{2(\alpha-1)}\psi(t)\int^t_0W(t-t')\partial_x[(\psi(t')w(t'))^3]dt'.\end{align}
Then, the proof of the theorem is summarized in four steps. Indeed, in the first two we prove that $\Phi_1$ is a contraction in the following ball of $X^{s,b}(\R\times\R)$:
\bel{boladem}\mathcal{B}\equiv\mathcal{B}(3c_0r_w):=\{w\in X^{s,b}:~~||w||_{X^{s,b}}\leq3c_0r_w\}.\eeq
1.- If $w\in\mathcal{B}$, $\alpha>2$ and combining (\ref{lema1}),~ (\ref{lema3}),~(\ref{bilinearkdv}),~(\ref{lematrilinear}), the following inequalities hold:
\begin{align}\nonumber
||\Phi_1(w)||_{X^{s,b}}&\leq c_0||w_0||_{H^{s}}+c_1\sigma\lambda^{\alpha-2}
||\psi(t)^2\partial_x(w^2(x,t))||_{X^{s,b-1}}\\
\nonumber &
+c_2\lambda^{2(\alpha-1)}||\psi(t)^3\partial_x(w^3(x,t))||_{X^{s,b-1}}\\
\nonumber &\leq c_0||w_0||_{H^{s}} + c\cdot
c_1\sigma\lambda^{\alpha-2}||w(x,t)||^2_{X^{s,b}}
+c\cdot c_2\lambda^{2(\alpha-1)}||w(x,t)||^3_{X^{s,b}}\\
\nonumber
&\leq c_0r_w + c\cdot c_1\sigma\lambda^{\alpha-2}(3c_0r_w)^2+c\cdot
c_2\lambda^{2(\alpha-1)}(3c_0r_w)^3\\
\nonumber
&\leq c_0r_w\left\{ 1 +
c_0\lambda^{\alpha-2}r_w+c_0^2\lambda^{2(\alpha-1)}r_w^2\right\}\\
\label{LWPfdem1}
&\leq 3c_0r_w,
\end{align}
where the last inequality is verified by choosing $\lambda$ which satisfy the following two conditions:
\begin{subequations}\label{condL}
 \begin{align}
   \label{condL01}r_w\lambda^{\alpha-2}c_0\leq1/4,\\
 \label{condL02}r_w^2\lambda^{2(\alpha-1)}c_0^2\leq1/4.
 \end{align}
\end{subequations}
Then, if we select $\lambda_0$ as the minimum value of \eqref{condL} and choose $\lambda$ as
\bel{Lmin}\lambda\leq\lambda_0=\min(d_1||w_0||^{\frac{-1}{\alpha-2}}_{H^s},d_2||w_0||^{\frac{-1 } {(\alpha-1)}} _{H^s}),~~d_1=(\frac{1}{4c_0})^{\frac{1}{\alpha-2}},d_2=(\frac{1}{4c_0^2})^{\frac{1}{2(\alpha-1)}},~~\alpha>2.\eeq
the conditions \eqref{condL} will be satisfied. 
In concluding,
\bel{contraccion}\Phi_1(\mathcal{B})\subseteq\mathcal{B}.\eeq
2.- With the same ideas as above (that is,  combining (\ref{lema1}),~(\ref{lema2}),~ (\ref{lema3}),
~ (\ref{bilinearkdv}),~(\ref{lematrilinear} and \eqref{LWPfdem1}), if $w,~\tilde{w}\in\mathcal{B}$ then we have 
%
%
%
\begin{align}\nonumber
||\Phi_1(w)-&\Phi_1(\tilde{w})||_{X^{s,b}}\leq c_0\sigma\lambda^{\alpha-2}
||\psi(t)^2\partial_x(w^2(x,t)-\tilde{w}^2(x,t))||_{X^{s,b-1}}\\
\nonumber
 &+c\lambda^{2(\alpha-1)}||\psi(t)^3\partial_x(w^3(x,t)-\tilde{w}^3(x,t))||_{X^{s,b-1}}\\
\nonumber
 &=  c_0\sigma\lambda^{\alpha-2}||\psi(t)^2\partial_x[(w(x,t)-\tilde{w}(x,t))(w(x,t)+\tilde{w}(x,t))]||_{X^{s,b-1}}\\
 \label{Phicontraccion}
&+c\lambda^{2(\alpha-1)}||\psi(t)^3\partial_x[(w(x,t)-\tilde{w}(x,t))(w(x,t)^2+\tilde{w}^2(x,t)+w(x,t)\tilde{w}(x,t))]||_{X^{s,b-1}}.
\end{align}

Now taking into account
\begin{itemize}
        \item [(i)]Bilinear estimate: given $u_1=w-\tilde{w},~~u_2=w+\tilde{w}$ , we have, using the bilinear estimate (see lemma \ref{bilinearkdv})
\begin{align}\nonumber||\psi(t)^2&\partial_x(w^2(x,t)-\tilde{w}^2(x,t))||_{X^{s,b-1}}=||\psi(t)^2\partial_x[u_1u_2]||_{X^{s,b-1}}\\
\nonumber& 
\leq
c_1||u_1||_{X^{s,b}}||u_2||_{X^{s,b}}
\leq c_1||w-\tilde{w}||_{X^{s,b}}\{||w||_{X^{s,b}}+||\tilde{w}||_{X^{s,b}}\}\\\label{Phicontraccion1}
&\leq 3\cdot2c_0c_1\cdot r_w||w-\tilde{w}||_{X^{s,b}}.\end{align}
\item [(ii)]Trilinear estimate: given $u_1=w-\tilde{w},~v_1=w,~v_2=\tilde{w},~u_2=w+\tilde{w}$, we have
 $$\begin{array}{llll}||\psi(t)^3\partial_x(w^3(x,t)-\tilde{w}^3(x,t))||_{X^{s,b-1}}\\\\
 
=||\psi(t)^3\partial_x[(w-\tilde{w})(w^2+\tilde{w}^2+w\tilde{w})]||_{X^{s,b-1}}\\\\

=||\psi(t)^3\partial_x[(w-\tilde{w})(w+\tilde{w})^2-(w-\tilde{w})w\tilde{w}]||_{X^{s,b-1}}\\\\

=||\psi(t)^3\partial_x[u_1u_2u_2]-\psi(t)^3\partial_x[u_1v_1v_2]||_{X^{s,b-1}}.\end{array}$$
and using the trilinear estimate (see lemma \ref{lematrilinear}), we obtain  
%




\begin{align}\nonumber||\psi(t)^3\partial_x(w^3(x,t)-\tilde{w}^3(x,t))||_{X^{s,b-1}}&
=||\psi(t)^3\partial_x[u_1u_2u_2]-\psi(t)^3\partial_x[u_1v_1v_2]||_{X^{s,b-1}}\\
\nonumber& 
\leq
c_2||u_1||_{X^{s,b}}||u_2||^2_{X^{s,b}} + c_2||u_1||_{X^{s,b}}||v_1||_{X^{s,b}}||v_2||_{X^{s,b}}\\
\nonumber
&\leq c_2||u_1||_{X^{s,b}}\{(||w||_{X^{s,b}}+||\tilde{w}||_{X^{s,b}})^2 + ||w||_{X^{s,b}}||\tilde{w}||_{X^{s,b}}\}\\\label{Phicontraccion2}
&\leq c_2||w-\tilde{w}||_{X^{s,b}}\{2(3c_0\cdot r_w)^2 + (3c_0\cdot r_w)^2\}\\
\nonumber&= 27c^2_0c_2\cdot r_w^2||w-\tilde{w}||_{X^{s,b}}.\end{align}
\end{itemize}
Hence \eqref{Phicontraccion} simplifies as 
\bel{Phicontraccion3}\begin{array}{llll}
||\Phi_1(w)-\Phi_1(\tilde{w})||_{X^{s,b}}\leq6c~c_0~c_1\cdot\sigma\lambda^{\alpha-2}r_w||w-\tilde{w}||_{X^{s,b}}\\\\

+ 27c~c^2_0~c_2\lambda^{2(\alpha-1)}  r_w^2||w-\tilde{w}||_{X^{s,b}}\\\\

\leq (c_0\lambda^{\alpha-2}r_w+c_0^2\lambda^{2(\alpha-1)} r_w^2)||w-\tilde{w}||_{X^{s,b}}\\\\
\leq \frac{1}{2} ||w-\tilde{w}||_{X^{s,b}},\end{array}\eeq
where the last inequality holds whenever we choose $\lambda$ less than the minimun value $\lambda_0$ in \eqref{Lmin}. Therefore, $\Phi_1$ is a contraction in $\mathcal{B}$ and then, there exists a unique $w\in\mathcal{B}(3c_0r_w)$ which satisfies
\bel{solucionPVI}
\begin{array}{ll} \psi(t)w(t)= \psi(t)\{W(t)w_0-6\sigma\lambda^{\alpha-2}\int^t_0W(t-t')\partial_x[(\psi(t')w(t'))^2]dt'\\\\
-2\lambda^{2(\alpha-1)}\int^t_0W(t-t')\partial_x[(\psi(t')w(t'))^3]dt'\},\end{array}\eeq
and in the temporal interval $[-1,1]$, $w(\cdot)$ is  solution of the IVP \eqref{GEw}. Recall that with the same argument as in \eqref{Phicontraccion3} it is also possible to prove that  
\bel{depdatoinicial}||w-\tilde{w}||_{X^{s,b}}\leq c~||w_0-\tilde{w}_0||_{H^{s}}.\eeq
3.- Before proving the property of persistence  \eqref{persistencia1}, we must show that the solution $w$ goes to the initial data $w_0$ in the classical $H^s$ norm, when $t\rightarrow0$. In this way, we must check that the integral terms on the right of the  following inequality go to $0$ when $t\rightarrow0$:
\bel{0persist}\begin{array}{l}||w(t)-\psi(t)W(t)w_0||_{H^{s}}          
\leq6\sigma\lambda^{\alpha-2}||\psi(t)\int^t_0W(t-t')\partial_x[(\psi(t')w(t'))^2] dt'||_{H^ {s}}\\\\
+2\lambda^{2(\alpha-1)}||\psi(t)\int^t_0W(t-t')\partial_x[ (\psi(t')w(t'))^3]dt'||_{ H^{s}},\end{array}\eeq
 To this end, we rewrite such terms in order to apply Lemma \ref{lemaa} \eqref{lema4} and analyze the behavior when $\eta\rightarrow0$ of the following general expression:
\bel{gPersist}||\psi(\eta)\int^\eta_0W(\eta-t')F(\cdot,t') dt'||_{H^ {s}}.\eeq
Taking  $t''=\frac{t'}{\eta}$, recalling  $t''\equiv t'$ and with the change (w.l.g.) of the cut off function from
$\psi(\eta)$ to $\psi(1)$, the expression \eqref{gPersist} is rewritten as follows:
 \bel{gPersist0}\begin{array}{l}||\psi(1)\int^1_0W(\eta(1-t'))F(\cdot,\eta t')
\eta dt'||_{H^ {s}}\\\\
=(\int|\psi(1)\int^1_0W(\eta(1-t'))(D^sF(\cdot,\eta
t'))(x)\eta dt'|^2dx)^{1/2}.\end{array} \eeq
With the change of variable $x=\eta^{1/3}y$, we rewrite \eqref{gPersist0} as

\bel{gPersist111}\begin{array}{l}(\int|\psi(1)\int^1_0W(\eta(1-t'))(D^sF(\cdot,
\eta t'))(\eta^{1/3} y)\eta dt'|^2\eta^{1/3}dy)^{ 1/2}\\\\
=\eta^{1+\frac{1}{6}}(\int|\psi(1)\int^1_0W(\eta(1-t'))(D^sF(\cdot,
\eta t'))(\eta^{1/3} y)dt'|^2dy)^{ 1/2}.\end{array}\eeq
Now, we analyze the integrand $W(\eta(1-t'))(D^sF(\cdot,\eta t'))(\eta^{1/3} y)$ of \eqref{gPersist111}. For that, we define
\bel{G1}G_1(y,\eta t')=D^sF(\eta^{1/3}y,\eta t')\equiv
D^sF_\eta (y,t'), \text{ where } F_\eta (y,t')=F(\eta^{1/3}y,\eta t').\eeq
%

%
In this way, $W(1-t')G_1(\cdot,\eta
t')(y)=W(\eta(1-t'))(D^sF(\cdot,\eta t'))(\eta^{1/3} y).$
Define
\bel{G2}G_2(y,t')=G_1(y,\eta t'),\eeq
\bel{G3}G_3(y,t')=D^{-s}G_2(y,t').\eeq
so that
\bel{Persist3}\begin{array}{l}D^sW(1-t')G_3(\cdot,t')=W(1-t')D^sG_3(\cdot,
t')=W(1-t')G_2(\cdot,t')\\\\
=W(1-t')G_1(\cdot,\eta t'),\end{array}\eeq
Now, we are able to apply the estimate Lemma \ref{lemaa} \eqref{lema4},
%
to obtain
\bel{gPersist1}\begin{array}{lll}\eta^{1+\frac{1}{6}}
(\int|\psi(1)\int^1_0W(\eta(1-t'))(D^sF(\cdot, \eta t'))(\eta^{1/3} y)dt'|^2dy)^{ 1/2}\\\\
=\eta^{1+\frac{1}{6}} \|\psi(1)\int^1_0W(1-t')G_3(\cdot,t')dt'\|_{H^s}\leq c\eta^{1+\frac{1}{6}}\|G_3\|_{X^{s,b-1}}\\\\
=c\eta^{1+\frac{1}{6}}\|D^{-s}G_2\|_{X^{s,b-1}}=c\eta^{1+\frac{1}{6}} \|G_2\|_{X^{0,b-1}}=c\eta^{1+\frac{1}{6}}\|D^sF_\eta \|_{X^{0,b-1}}\\\\
=c\eta^{1+\frac{1}{6}}\|F_\eta \|_{X^{s,b-1}}.\end{array}\eeq
Therefore, we must calculate
$\|F_\eta\|_{X^{s,b-1}}=\|F(\eta^{1/3}y,\eta t')\|_{X^{s,b-1}}$.\\
\bel{normaXTras}\begin{array}{lll}\|F(\eta^{1/3}y,\eta t')\|_{X^{s,b-1}}
=(\int_{\R}\int_{\R}(1+|\tau-\xi³|)^{2(b-1)}(1+|\xi|)^{2s}|\widehat{F(\eta^{1/3}y, \eta t')}(\xi,\tau))|^2d\xi d\tau)^{1/2}.\end{array}\eeq
Since
$\widehat{F(\eta^{1/3}y,\eta t')}(\xi,\tau)=\eta^{-4/3}\hat{F}(\eta^{-1/3}\xi,\eta^{-1} \tau)$,

\bel{normaXTras1}\begin{array}{l}\|F(\eta^{1/3}y,\eta t')\|_{X^{s,b-1}}=(\int_{\R}\int_{\R}(1+|\tau-\xi³|)^{2(b-1)}
(1+|\xi|)^ { 2s} |\eta^{-4/3}\hat{F}(\eta^{-1/3}\xi,\eta^{-1}\tau)|^2d\xi
d\tau)^{1/2}.\end{array} \eeq
With the change of variables,~$\xi'=\eta^{-1/3}\xi,~~\tau'=\eta^{-1}\tau$, the above simplifies as follows
\bel{normaXTras2}\begin{array}{l}\|F(\eta^{1/3}y,\eta t')\|_{X^{s,b-1}}=(\int_{\R}\int_{\R}(1+\eta|\tau'-\xi'³|)^{
2(b-1)}(1+\eta^{1/3}|\xi'|)^ { 2s} |\hat{F}(\xi',\tau')|^2 
\eta^{-8/3}\eta^{4/3}d\xi'd\tau')^{1/2}\\\\
=\eta^{-2/3}(\int_{\R}\int_{\R}(1+\eta|\tau'-\xi'³|)^{
2(b-1)}(1+\eta^{1/3}|\xi'|)^ { 2s} |\hat{F}(\xi',\tau')|^2 
d\xi'd\tau')^{1/2}.\end{array} \eeq
Now, we must consider the following different cases,

\begin{enumerate}
 \item[(1)]
$\begin{cases}\eta^{1/3}|\xi'|\leq1,\\\eta|\tau'-\xi'³|\leq1.\end{cases}$
In this case, $\|F_\eta\|_{X^{s,b-1}}\leq \eta^{-2/3}\|F\|_{X^{s,b-1}}.$\\

 \item[(2)]
$\begin{cases}\eta^{1/3}|\xi'|\leq1,\\\eta|\tau'-\xi'³|\geq1.\end{cases}$
In this case, $\|F_\eta\|_{X^{s,b-1}}\leq \eta^{b-1-2/3}\|F\|_{X^{s,b-1}}.$\\

\item[(3)]$\begin{cases}\eta^{1/3}|\xi'|\geq1,\\\eta|\tau'-\xi'³|
\leq1.\end{cases}$
In this case, $\|F_\eta \|_{X^{s,b-1}}\leq \eta^{s/3-2/3}\|F\|_{X^{s,b-1}}.$\\

\item[(4)]$\begin{cases}\eta^{1/3}|\xi'|\geq1,\\\eta|\tau'-\xi'³|          
\geq1.\end{cases}$
In this case, $\|F_\eta\|_{X^{s,b-1}}\leq \eta^{s/3+b-1-2/3}\|F\|_{X^{s,b-1}}.$\\
\end{enumerate}
Recall that $1/2<b<1$.  We use the more restrictive upper bound which corresponds to the second case to conclude that the estimate \eqref{gPersist1} remains as follows:

\bel{gPersist2}\begin{array}{l}c\eta^{1+1/6}\|F_\eta\|_{X^{s,b-1}}\leq
c\eta^{1+1/6+b-1-2/3}\|F\|_{X^{s,b-1}}\\\\
=c\eta^{b-1/2}\|F\|_{X^{s,b-1}}.\end{array}\eeq
Then, \eqref{gPersist} goes to $0$ when $\eta \rightarrow0$ and, using the continuity of the free operator $W(t)$,

\bel{contiDatoini}\lim_{t\rightarrow0}||w(t)-\psi(t)W(t)w_0||_{H^s}=0.\eeq\\
4.- We now prove the persistence property \eqref{persistencia1}, i.e. $w\in C([-1,1]:H^s(\mathbb{R}))$. First, we need to prove the continuity of the norm of the $X^{s,b}$ space but this follows a similar argument as in the proof of the continuity of the initial data in the $H^s$ norm \ref{0persist}, therefore it will be omitted. From this, the persistence property in $H^s$ is a direct consequence of the embedding  $X^{s,b} \subset C_tH^s ~~\text{for}~~ b > 1/2$ (see \cite[p.156, corollary 7.3]{LiPo}).\begin{flushright}$\Box$\end{flushright}
%
{\bf Proof of Theorem \ref{LWPf}.}\\
%
From the unicity and local existence, at time $t=T>0$ of the IVP \eqref{GEw} for  $w$, we get  unicity and  local existence, at time $T=\lambda^{-3}$, of the IVP \eqref{GE0} for $v$, whenever we determine for which values of $\lambda$, depending on  $||v_0||_{H^s}$, conditions \eqref{condL01} and  \eqref{condL02} are verified. In this way, we must calculate the norm $r_w=||w_0||_{H^{s}}$, taking now into account  that $w_0(x)=\lambda^{-\alpha}v_0(\lambda^{-1}x)$, using \eqref{escala1}. By definition
\bel{normaW1}\begin{array}{l}
   r_w=||w_0||_{H^{s}}=(
\int_{\R}(1+|\xi|)^{2s}|\widehat{w_0(x)}(\xi)|^2d\xi)^{1/2},
  \end{array}\eeq
where,
%
\begin{align}\nonumber \widehat{w_0(x)}(\xi)&=\int_{\R}e^{-ix\xi}w_0(x)dx=\int_{\R}e^{-ix\xi}\lambda^
{-\alpha}v_0(\lambda^{-1}x)dx\\
\nonumber&=\int_{\R}e^{-i(\lambda^{-1}x)\lambda\xi}
\lambda^{-\alpha}v_0(\lambda^{-1}x)\lambda
d(\lambda^{-1}x)=\lambda^{1-\alpha}\int_{\R}e^{-i(\lambda^{-1}x)\lambda\xi}
v_0(\lambda^{-1}x)d(\lambda^{-1}x)\\
\nonumber&=\lambda^{1-\alpha}\hat{v}_0(\lambda\xi).\end{align}
Then, \eqref{normaW1} is rewritten as follows

%
\begin{align}\nonumber  r_w&=(\int_{\R}(1+|\xi|)^{2s}|\widehat{w_0(x)}(\xi)|^2d\xi)^{1/2}=
(\int_{\R}(1+|\xi|)^{2s}|\lambda^{1-\alpha}\hat{v}_0(\lambda\xi)|^2d\xi)^{1/2}\\
\nonumber&=(\lambda^{1-\alpha}\int_{\R}(1+\lambda^{-1}|\lambda\xi|)^{2s}|\hat{v}
_0(\lambda\xi)|^2\lambda^{-1}d(\lambda\xi))^ { 1/2 }\\
\nonumber&=\lambda^{1-1/2-\alpha}(\int_{\R}(1+\lambda^{-1}|\lambda\xi|)^{2s}|\hat{v}
_0(\lambda\xi)|^2d(\lambda\xi))^{ 1/2 }\\\label{normaW2}
&=\lambda^{1/2-\alpha}(\int_{\R}(1+\lambda^{-1}|\lambda\xi|)^{2s}|\hat{v}
_0(\lambda\xi)|^2d(\lambda\xi))^{ 1/2 }.\end{align}

Now denoting  $r_v=||v_0||_{H^{s}}$,  the norm $r_w=||w_0||_{H^{s}}$ can be estimated in terms of $r_v$, splitting first the  last integral in \eqref{normaW2}  and  estimating next:

\bel{normaW3}   
\begin{array}{l}
(a)~~\lambda^{1/2-\alpha}(\int_{\lambda^{-1}|\eta|\leq1}|v_0(\eta)|^2d\eta)^{
1/2 }\leq \lambda^{1/2-\alpha}||v_0||_{L^2(\R)}\leq\lambda^{1/2-\alpha}r_v,\\\\
(b)~~\lambda^{1/2-\alpha-s}(\int_{\lambda^{-1}|\eta|\geq1}
|\eta|^{2s}|v_0(\eta)|^2d\eta)^
{1/2}\leq\lambda^{1/2-\alpha-s}||v_0||_{\dot{H}^s(\R)}\leq\lambda^{1/2-\alpha-s}
r_v. \end{array}\eeq
Now, taking into account  the latter bounds for the norm  $r_w$ when $\lambda^{-1}|\eta|\leq1$ and $\lambda^{-1}|\eta|\geq1$ respectively and recalling that $s>1/4$, we determine for which values of $\lambda$ the conditions \eqref{condL} are satisfied. 

[(i)] Condition \eqref{condL01}: $r_w\lambda^{\alpha-2}c_0\leq1/4.$\\
\begin{itemize}
 \item
[(i.a)]$$\lambda^{1/2-\alpha}r_v\lambda^{\alpha-2 }c_0\leq1/4\rightarrow
\lambda^{-3/2}\leq\frac{1}{4c_0r_v}.$$
\item
[(i.b)]$$\lambda^{1/2-\alpha-s}r_v\lambda^{\alpha-2}c_0 \leq1/4\rightarrow
\lambda^{-3/2-s}\leq\frac{1}{4c_0r_v}.$$
\end{itemize}

[(ii)] Condition \eqref{condL02}: $r^2_w\lambda^{2(\alpha-1)}c_0^2\leq1/4.$\\
\begin{itemize}
 \item
[(ii.a)]$$\lambda^{1-2\alpha}r^2_v\lambda^{2(\alpha-1)}c_0^2 \leq1/4\rightarrow
\lambda^{-1}\leq\frac{1}{4c_0^2(r_v)^2}.$$
\item
[(ii.b)]$$\lambda^{1-2\alpha-2s}r^2_v\lambda^{2(\alpha-1)}
c_0^2\leq1/4\rightarrow\lambda^{-1-2s}\leq\frac{1}{4c_0^2(r_v)^2}.$$
\end{itemize}
In this way, we distinguish two different cases:

\begin{itemize}
 \item[1.] $$\frac{1}{4c_0r_v}\leq1$$
From [(i)], if (i.a) is satisfied, (i.b)  will also be satisfied. Then,
$$\lambda^{-3/2}\leq\frac{1}{4c_0r_v}.$$\\
From [(ii)], if (ii.a) is satisfied, (ii.b) will also be satisfied. So that
$$\lambda^{-1}\leq\frac{1}{4c_0^2(r_v)^2}\Rightarrow\lambda^{-3/2}\leq\frac{1}{
2c_0^2r_v^2(4c_0r_v) }$$\\
Comparing these two estimates on $\lambda$, the most restrictive is

\bel{condL11}\lambda^{-1}\leq\frac{1}{4c_0^2(r_v)^2}.\eeq
\item[2.] $$\frac{1}{4c_0^2(r_v)^2}\geq1$$
From [(i)], if (i.b) is satisfied, (i.a) will also be satisfied. Therefore, 
$$ \lambda^{-3/2-s}\leq\frac{1}{4c_0r_v}.$$\\
From [(ii)], if (ii.b) is satisfied, (ii.a) will also be satisfied. Then,

$$\lambda^{-1-2s}\leq\frac{1}{4c_0^2(r_v)^2}.$$\\
From these two estimates on $\lambda$, the most restrictive is

\bel{condL12}\lambda^{-3/2-s}\leq\frac{1}{4c_0r_v}.\eeq
\end{itemize}

Hence, whenever $\lambda$ satisfies \eqref{condL11} or \eqref{condL12}, the application $\Phi_{T=\lambda^{-3}}$ is contractive in $\mathcal{B}_v$, the  ball associated to the ball $\mathcal{B}$ given in \eqref{boladem} through the scaling relation. So, there exists  a  unique $v\in\mathcal{B}_v$ solution of the IVP \eqref{GE0}.\\

Note that the local existence time is explicit  from the relation $\lambda³ T=1$ and from the values of $\lambda$ \eqref{condL11} and \eqref{condL12}.
That is,
\bel{T1}T=\lambda^{-3}\leq(\frac{1}{4c_0^2||v_0||^2_{H^s}})^{3},~~||v_0||_{H^s}
\geq1,\eeq
\bel{T2}T=\lambda^{-3}\leq(\frac{1}{4c_0||v_0||_{H^s}})^{\frac{6}{3+2s}},
~~||v_0||_ { H^s}\leq1.\eeq
With the same argument used in \eqref{Phicontraccion3}, it verifies that

\bel{depdatoinicial1}||v-\tilde{v}||_{X^{s,b}}\leq
c~||v_0-\tilde{v}_0||_{H^{s}}.\eeq
The property of persistence \eqref{persistencia1} for $v$, i.e.

$$v\in C([-T,T]:H^s(\mathbb{R})),$$
follows directly from the  property of persistence for the auxiliar function $w$. The proof is now complete. \begin{flushright}$\Box$\end{flushright}

\subsection{\textbf{Global Theory}}\label{globalWP}
The global well-posedness for the IVP of the GE \eqref{GE0} with initial data in $H^1(\R)$ is stated in the following theorem:

\begin{thm}\label{GWPdef2} Let be $u_0\in H^{1}(\R)$ and  $u$ the corresponding local solution for the IVP of the GE \eqref{GE0} given by the theorem  \ref{LWPf}. Then we extend the solution for all $t>0$, that is

\bel{globaltime2}u\in C(\R:H^{1}(\R)).\eeq
\end{thm}
{\it Proof.} It follows using the local existence result and standard techniques of Gagliardo-Nirenberg inequalities and conservation of the energy.\begin{flushright}$\Box$\end{flushright}

\section{\textbf{Stability of solitons}}\label{Stability}

We consider here the question of the stability of solitons of the Gardner equation \eqref{GE0} under small perturbations in $H^1(\R)$. Since the mKdV and  Gardner equations are closely related, this question is equivalent to study the stability under small perturbations in  $H^1(\R)$ of solitons with non bounded mean of the focusing and defocusing mKdV (positive/negative nonlinearity, respectively), 


\bel{mkdvfdef}\frac{\partial }{\partial t}k(x,t) + \frac{\partial}{\partial x^3}k(x,t) \pm 2\frac{\partial}{\partial x}(k^3(x,t))=0.\eeq

In this section, we compute necessary conditions, given by P. Zhidkov (see \cite{Z}) in a general framework,  to obtain the stability of solitons. In such work, P. Zhidkov states a general theorem about the stability of solitons of the gKdV equation in $H^2(\R)$, vanishing in the boundary. In our case, the stability result is centered on the stability of solitons of the Gardner equation.

        
        \subsection{Existence of travelling wave solutions}\sectionmark{Existence of travelling wave solutions}\label{existence}
        We look for solutions of the focusing mKdV 
        \bel{mkdvf}\frac{\partial }{\partial t}k(x,t) + \frac{\partial}{\partial x^3}k(x,t) + 2\frac{\partial}{\partial x}(k^3(x,t))=0,\eeq     
        of the following type
        \bel{ansatzf2}k(x,t) = \sigma + \phi(x-c_\sigma t),~~c_\sigma>0,~\phi(\pm\infty)=0.\eeq
Introducing this ansatz in \eqref{mkdvf}, we obtain   $-c_\sigma\phi'+\phi^{'''} + 2((\sigma+\phi)^3)'=0$. Integrating and taking into account that $\phi(\pm\infty)=0$, we arrive to
        
        \bel{ivpf}\begin{cases}\phi^{''} + 2\phi^3 +6\sigma\phi^2 + (6\sigma^2-c_\sigma)\phi=0\\\phi(\pm\infty)=0.\end{cases}\eeq
Now multipliying \eqref{ivpf} by the  integrating factor $\phi'$,
 we arrive to
        
\begin{figure}[!htb]
\centering
\includegraphics[width=12cm,height=4.5cm]{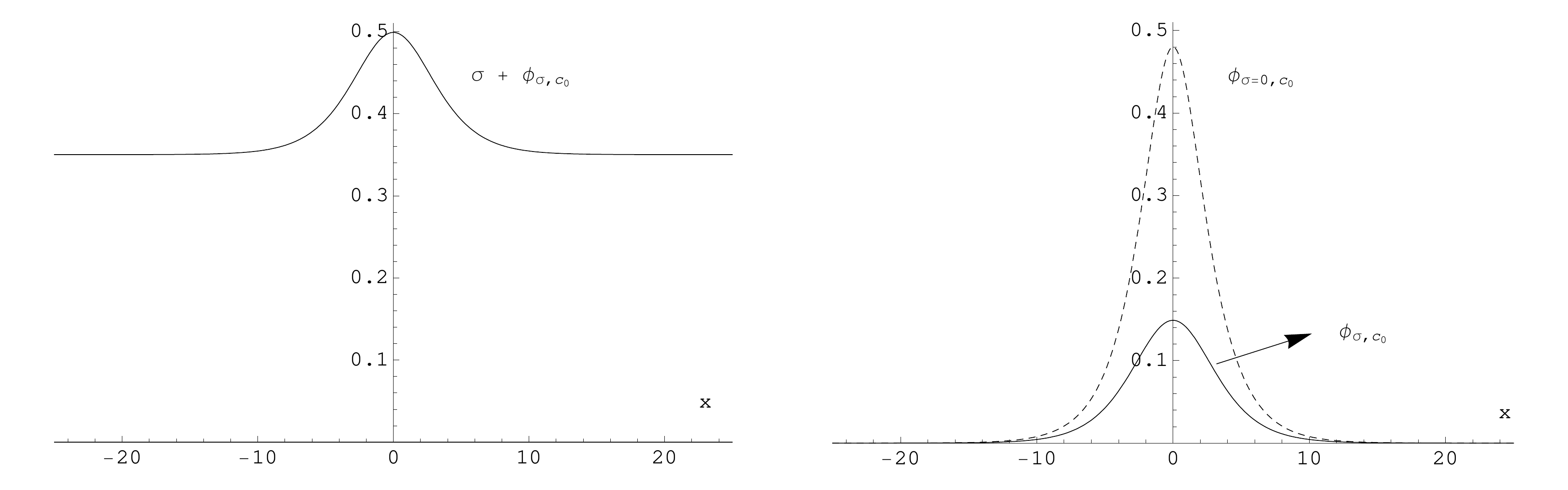}
\small{\caption{Left. Profile of $\sigma+\phi_{\sigma,c_0}$ with $\sigma=0.35,~~c_0=0.23$. Right.
Comparison between humps of the  soliton of mKdV (dashed line, $\phi_{\sigma=0,c_0}$) and the soliton of GE  (solid line, $\phi_{\sigma,c_0}$) with $\sigma=0.35$, and both with $c_0=0.23$.\label{solitonGE}}}
\end{figure}

        \bel{edof}(\phi')^2+\phi^4+4\sigma\phi^3+(6\sigma^2-c_\sigma)\phi^2=0.\eeq
This ODE can be solved explicitly and we obtain (see figure \ref{solitonGE})
\bel{fexplicita2}\phi_{\sigma,c_0}(x-c_\sigma(\sigma,c_0)t)=\frac{c_0}{2\sigma+\sqrt{4\sigma^2+c_0}\cosh(\sqrt{c_0}(x-(6\sigma^2+c_0)t))}.\eeq
        $$\begin{array}{ll}c_\sigma(\sigma,c_0)=6\sigma^2+c_0.\end{array}$$     

We can also characterize the soliton of the Gardner equation as the minimum of the following Lyapunov functional:
if we denote    
        \bel{conservada12}\begin{array}{ll}{\mathcal{E}}(f)=\int_{\mathbb{R}}\left\{f_x^2 - f^4 - 4\sigma f^3 \right\}dx,\\\\

 {\mathcal{F}}(f)=\frac{1}{2}\int_{\mathbb{R}}f^2dx,\end{array}\eeq
as the energy and the $L^2$ norm of the soliton of the Gardner equation \eqref{GEw}, perturbing the Lyapunov functional \eqref{funcionalf} 
        
        \bel{funcionalf}E(\xi)={\mathcal{E}}(\xi) + 2(c_\sigma-6\sigma^2){\mathcal{F}}(\xi) = \int_{\mathbb{R}}\left\{\xi_x^2 - \xi^4 - 4\sigma\xi^3 + (c_\sigma-6\sigma^2)\xi^2\right\}dx,\eeq
around its critical point, we obtain the linearized operator $L$ and the ODE satisfied by the soliton of the Gardner equation,  

        \bel{operadorf}L = -\partial_{xx} + c_\sigma-6(\sigma+\phi_{\sigma,c_0})^2.\eeq
        \bel{eqRf}-\phi_{\sigma,c_0}^{''}-2\phi_{\sigma,c_0}^3-6\sigma\phi_{\sigma,c_0}^2+(c_\sigma-6\sigma^2)\phi_{\sigma,c_0}=0.\eeq
        
The main results of this section are the following:\\
        \begin{thm}[focusing case]\label{estabilidadf} Let $\sigma\in\mathbb{R},~c_0\in(0,\infty)$ and $u_{\sigma,c_0}(x,t)=\sigma+\phi_{\sigma,c_0}(x-c_\sigma(\sigma,c_0)t)\in\dot{H}^1(\R)$, $c_\sigma(\sigma,c_0)=c_0+6\sigma^2,$ a solution of the focusing mKdV \eqref{mkdvf}, where $\phi_{\sigma,c_0}\in H^1(\R)$ satisfies (\ref{ivpf}). Then \\

$\forall\epsilon>0,~\exists\delta\equiv\delta(\epsilon,\sigma,c_0)>0$ and a  function $C^2(\mathbb{R})$,~~ $r:\mathbb{R}\rightarrow\mathbb{R}$, such that\\

if  $||u_0-(\sigma+\phi_{\sigma,c_0})||_{H^1(\R)}<\delta$, then\\ 
$$\sup_{t>0}\left\|u(\cdot,t)-u_{\sigma,c_0}(\cdot+r(t))\right\|_{H^{1}(\mathbb{R})}<\epsilon,$$
where $u(x,t)$  is the unique solution of the focusing mKdV equation with initial data $u_0=u(x,0)\in\dot{H}^1(\R)$ and where $\sup_{t}|r'(t)+(c_0+6\sigma^2)|\leq K\epsilon,~~K>0$.
\end{thm}

\begin{thm}[defocusing case]\label{estabilidaddef}
Let $\sigma\in\mathbb{R},~c_0\in(0,4\sigma^2)$ and $u_{\sigma,c_0}(x,t)=\sigma-\varphi_{\sigma,c_0}(x+c_\sigma(\sigma,c_0)t)\in\dot{H}^1(\R)$,~~$c_\sigma(\sigma,c_0)=6\sigma^2-c_0,$ a solution of the defocusing mKdV, where $\varphi_{\sigma,c_0}\in H^1(\R)$ is solution of the defocusing Gardner equation. Then\\ 
$\forall\epsilon>0,~\exists\delta\equiv\delta(\epsilon,\sigma,c_0)>0$ and a  function $C^2(\mathbb{R})$,~~ $r:\mathbb{R}\rightarrow\mathbb{R}$, such that\\

if $||u_0-(\sigma-\varphi_{\sigma,c_0})||_{H^1(\R)}<\delta$, then\\

$$\sup_{t>0}\left\|u(\cdot,t)-u_{\sigma,c_0}(\cdot+r(t))\right\|_{H^{1}(\mathbb{R})}<\epsilon,$$
where $u(x,t)$  is the unique solution of the focusing mKdV equation with initial data $u_0=u(x,0)\in\dot{H}^1(\R)$ and where $\sup_{t}|r'(t)+(c_0-6\sigma^2)|\leq K\epsilon,~~K>0$.
\end{thm}
Before explaining the main ideas behind the proof of this result, some remarks are in order.\\
        {\bf Remarks.}\\ 
        
\begin{enumerate}
        \item[(i)] We will fix the parameter $\sigma$ and we will choose the parameter $c_0$ in the interval $c_0\in(0,\infty)$  since we are only interested in real regular solutions. In this way, the aplication 
\bel{Phifoc}c_0\in(0,\infty)\stackrel{\Phi}{\longrightarrow}\phi_{\sigma,c_0}\in H^{1}(\mathbb{R})\eeq
         is $C^1(\mathbb{R^{+}}:H^{1}(\mathbb{R}))$.
        \item[(ii)] The soliton solution of the defocusing mKdV equation is easily obtained from the explicit expression \eqref{fexplicita2} in the focusing case, with the change $\sigma$ to $i\sigma$.\\
\end{enumerate}

{\bf Proof of Theorem \ref{estabilidadf}}. We use the standard techniques given by M.Weinstein \cite{We}, P.Zhidkov \cite{Z} and J.Angulo \cite{An1,An2} and then we do not give the details here. We only comment two main results in the proof:\\

\begin{enumerate}
        \item Convexity of the function $d(c_0)$: $d''(c_0)>0,~~\forall c_0\in(0,\infty)$, where
        \bel{convexi}d(c_0)=\int_{\mathbb{R}}\left\{\partial^2_x\phi_{\sigma,c_0} - \phi_{\sigma,c_0}^4 - 4\sigma\phi_{\sigma,c_0}^3 + (c_\sigma-6\sigma^2)\phi_{\sigma,c_0}^2\right\}dx.\eeq
        
        Since the Gardner equation is not invariant under scaling transformations, we can not proceed as usual to prove the convexity of the function \eqref{convexi}. But  we know explicitly the soliton solution $\phi_{\sigma,c_0}$ of the Gardner equation (see \eqref{fexplicita2}) and then by its relative simplicity, we can integrate directly (maybe with the help of a handbook of integrals, e.g. \cite{GrR}) obtaining
        
        $$d''(c_0)=2\int_{\mathbb{R}}\phi_{\sigma,c_0}\partial_{c_0}\phi_{\sigma,c_0}ds=\frac{\sqrt{c_0}}{4\sigma^2+c_0}>0,~~\forall c_0\in(0,\infty).$$

        \item The phase $r(t)$ and its velocity  $r'(t)$.\\
Define the function $F:\mathbb{R}^2\rightarrow\mathbb{R}$ given by,\\
\bel{funcionFf}\begin{array}{ll}F(r,t)=\frac{1}{2}\int_{\mathbb{R}}\left\{u(x,t)-(\sigma+\phi_{\sigma,c_0}(x+r))\right\}^2dx\\\\
=\frac{1}{2}\int_{\mathbb{R}}\left\{\psi(x,t)-\phi_{\sigma,c_0}(x+r)\right\}^2dx.\end{array}\eeq
Since $\phi_{\sigma,c_0}\in H^{\infty}(\R)$, $F$  is a  $C^{\infty}$ function in the $r$ variable and $C^1$ in the $t$ variable. With these properties for $F$ we define the following associated function,\\

\bel{funcionGf}
\begin{array}{ll}
G(r,t)=\frac{\partial F}{\partial r}=-\int_{\mathbb{R}}\left\{u(x,t)-(\sigma+\phi_{\sigma,c_0}(x+r))\right\}\partial_x\phi_{\sigma,c_0}(x+r)dx\\\\
=-\int_{\mathbb{R}}u(x,t)\partial_x\phi_{\sigma,c_0}(x+r)dx.\end{array}\eeq
It is easy to see that $G(0,0)=0$ and $\frac{\partial G}{\partial r}|_{(0,0)}\geq0$. Then, applying the implicit function theorem to the function $G$, there exist $T>0$ and a  $C^2((-T,T):\R), \text{ function } r:(-T,T)\rightarrow\R,~~r(0)=0$, such that (from $G$ and its time derivatives, we check that the function $r(t)$ is $C^2$)
        
        \bel{consecTFimplicta}G(r(t),t)=0,~~\forall t\in(-T,T).\eeq
        We choose $T$ in a maximal way. Then derivating implicitly \eqref{consecTFimplicta}, we get
        
        $$\frac{d}{dt}G(r(t),t)=0=\frac{\partial G}{\partial r}r'(t) + \frac{\partial G}{\partial t}\Rightarrow$$
        
        \bel{obtencionR}r'(t)=-\frac{\frac{\partial G}{\partial t}}{\frac{\partial G}{\partial r}}.\eeq
        Then, once we calculate $\frac{\partial G}{\partial t},~~\frac{\partial G}{\partial r}$, and substituting in \eqref{obtencionR}, we obtain the ODE satisfied by the phase $r(t)$

\bel{ODEr'f}\begin{cases}r'(t)=-(c_0+6\sigma^2)-\frac{\int_{\mathbb{R}}h\{-12v(v')^2+6hvv^{''}+2h^2v^{''}\}dx}{\int_{\mathbb{R}}\{-(v^{'})^2+hv^{''}\}},\\\\
r(0)=0,\end{cases}\eeq\\
where, $h(x,t)=u(x,t)-(\sigma+\phi_{\sigma,c_0}(x+r(t)))~~~\wedge~~~v(x,t)=\sigma+\phi_{\sigma,c_0}(x+r(t)).$
\begin{flushright}$\Box$\end{flushright}
\end{enumerate}
{\bf Proof of Theorem \ref{estabilidaddef}}. It follows the same steps than in the focusing case (up to the obvious change in the linearized operator and the speed of the soliton), with the important difference of the existence interval for the parameter $c_0\in(0,4\sigma^2)$.\begin{flushright}$\Box$\end{flushright}

\bigskip

\noindent
{\bf Acknowdlegments.} The author wishes to express his sincere thanks to Professor Luis Vega  for his continuous encouragement and useful discussions during the elaboration of this work.

\bigskip



\end{document}